\documentclass[12pt,letterpaper,fleqn]{article}
\usepackage[nosumlimits]{amsmath}
\usepackage{amssymb}
\usepackage{bm}
\usepackage{color}
\usepackage{natbib}
\usepackage[dvipdfmx]{graphicx}
\usepackage{lscape}
\usepackage{threeparttable}
\usepackage{here}
\usepackage{amsmath,amssymb}
\usepackage[top=30truemm,bottom=30truemm,left=25truemm,right=25truemm]{geometry}
\renewcommand{\baselinestretch}{1.0}
\begin{document}
\setlength{\footnotesep}{0.6cm}
\begin{titlepage}
\renewcommand{\baselinestretch}{1.0}
\title{Positive definiteness of the asymptotic covariance matrix of  
OLS estimators in parsimonious regressions%; comments on  
%``Testing a large set of zero restrictions in regression models, with an application to mixed frequency Granger causality'' by \citet{GHM01} 
\footnote{The author would like to thank 
Kaiji Motegi for his helpful comments. 
Of course, any remaining errors are the author's own.}
}
\author{ Daisuke Nagakura  \\ \\
Faculty of Economics, Keio University,  }
   \date{ This version: October 26, 2020
\\\
First version: October 13, 2020
}
\maketitle

%\newpage
\begin{center}
\textbf{Abstract}
\end{center}
Recently, \citet{GHM01} proposed a test for examining 
whether a large number of coefficients in linear regression models 
are all zero. The test is called the ``max test.'' 
The test statistic is calculated by first 
running multiple ordinary least squares (OLS) regressions, each including  
only one of key regressors, 
whose coefficients are supposed to be zero under the null, and then 
taking the maximum value of the squared OLS coefficient estimates 
of those key regressors. 
They called these regressions ``parsimonious regressions.''  
This paper answers a question raised in their Remark 2.4;  
whether the asymptotic covariance matrix of the OLS estimators 
in the parsimonious regressions is generally 
positive definite. The paper shows that it is generally positive definite, and  
the result may be utilized to 
facilitate the calculation of the simulated $p$ value necessary 
for implementing the max test.  
$ $\\ \\
Keywords: Max test; Linear regression; Positive definiteness 
Many covariates ; Parsimonious regression.
\\
JEL Classification: C12, C22, C51.
\end{titlepage}
\setcounter{page}{2}
\section{Introduction}
\newtheorem{lemma}{Lemma}
\newtheorem{theorem}{Theorem}
\newtheorem{corollary}{Corollary}
\newtheorem{proposition}{Proposition}
Several studies have considered the problem 
of testing a large number of restrictions in linear regression models 
when the number of regressors is large relative to the number of 
the observations.  
For instance, \citet{CAL01} showed that the usual \textit{F} 
test has an over-rejection 
tendency, and thus, proposed a modified \textit{F} test 
that has correct sizes even in this case. 
Similarly, \citet{ANT01} derived modified \textit{F}, 
likelihood ratio, and Lagrange multiplier tests 
under assumptions similar but different from those of \citet{CAL01} (see \citet{RIC01} for the bootstrapped versions of these tests). 
These tests are modifications of the well-known classical tests 
and are derived in the context of cross-sectional regressions with 
homoscedastic errors. 
Please refer to \citet{CJN01}, \citet{CJN02}, \citet{CKP01}, 
and \citet{CJM01} for related analyses. 

Recently, \citet{GHM01} (hereafter, GHM) proposed a new test for examining 
a large set of zero restrictions in linear regression models 
that possibly include time series variables and heteroscedastic errors. 
They called this new test the ``max test.'' 
Specifically, they considered the following regression model: 
\begin{equation}
\label{full}
y_t = \bm{z}_t^{\prime} \bm{a} + 
\bm{x}_t^{\prime} \bm{b} + \epsilon_t,\quad t=1,...,n,
\end{equation}
where $\bm{z}_t=[z_{1t}, ... , z_{pt}]^{\prime}$ and $\bm{x}_t=[x_{1t}, ... , x_{ht}]^{\prime}$ are 
vectors of explanatory variables with dimensions 
$p$ and $h$, respectively, $\bm{a}=[a_{1}, ... , a_{p}]^{\prime}$, 
$\bm{b}=[b_1, ... , b_h]^{\prime}$ are vectors of 
coefficients, and $\epsilon_t$, $t=1,...,n$, are error terms with 
zero mean and a finite second moment. 
Dimension $p$ is assumed to be small, while 
$h$ is large but finite. %The inequality $p+h \leq n$ is also assumed. 
It is well known that the usual Wald test      
for the null hypothesis $H_0: \bm{b}=\mathbf{0}_{h\times1}$ and against 
alternative hypothesis $H_1:  \bm{b} \not=\mathbf{0}_{h\times1}$  
has a severe over-rejection tendency in small sample due to parameter proliferation,  
where $\mathbf{0}_{a\times b}$ 
denotes the $a\times b$ vector with all elements being zero. 
Moreover, bootstrapping the Wald test corrects the finite sample size,  
reducing the power significantly.  
To overcome these difficulties, GHM proposed the max test, and  
their simulation experiments confirmed that the max test 
does not show the over rejection tendency of the usual Wald test,  
and is more powerful than the bootstrapped Wald test. 

To construct the max test statistic, we first run the following regressions for $i=1,..,h$, 
which GHM called ``parsimonious regressions'':
\begin{equation}
\label{PR}
y_t = \bm{z}_t^{\prime} \bm{\alpha}_i + 
x_{it} \beta_i + u_{it}, 
\end{equation}
where $u_{it}$ is an error term that depends on $i$. 
Let $\widehat{\beta}_{ni}$ be the OLS estimate of $\beta_i$ from the 
$i$th parsimonious regression in (\ref{PR}). GHM defined the max test statistic, $\widehat{\mathcal{T}}_n$, as:
\begin{equation}
\label{Tmax}
 \widehat{\mathcal{T}}_n=\max\{ (\sqrt{n}\widehat{\beta}_{n1})^2, ... , 
(\sqrt{n}\widehat{\beta}_{nh})^2\}.
\end{equation}
When $\widehat{\mathcal{T}}_n$ is sufficiently large, the null hypothesis is rejected. 

Let $\bm{V}$ denote the asymptotic covariance matrix of 
$\sqrt{n}\widehat{\bm{\beta}}_n$, where  
$\widehat{\bm{\beta}}_n = [\widehat{\beta}_{n1}, ...,\widehat{\beta}_{nh}]^{\prime}$. 
The asymptotic distribution of the max test statistic depends on $\bm{V}$,  
and is thus not pivotal.  
To implement the test, GHM utilized simulated $p$ value  
for $\widehat{\mathcal{T}}_n$, proposing to compute the simulated $p$ value 
in the following way. 
First, generate a sample $\{ \mathcal{N}_i^{(j)} \}_{i=1}^h$   
from the asymptotic distribution of $\sqrt{n}\widehat{\bm{\beta}}_n$ 
with estimated $\bm{V}$, or 
$N(\mathbf{0}_{h\times 1}, \widehat{\bm{V}}_n)$, where $\widehat{\bm{V}}_n$ is 
an estimate of $\bm{V}$, then, calculate 
$\widehat{\mathcal{T}}_n^{(j)}\equiv\max\{ \{ \mathcal{N}_i^{(j)2}\}_{i=1}^h \}$, 
repeat the procedure $M$ times to obtain  
$\widehat{\mathcal{T}}_n^{(j)}$, $j=1,...,M$, and finally compute the 
ratio of the number of $\widehat{\mathcal{T}}_n^{(j)}s$ that 
are greater than $\widehat{\mathcal{T}}_n$ to $M$ 
as the simulated $p$ value (see GHM for the computational details).
GHM proved that this simulated $p$ value converges to the 
true $p$ value and the test has the correct size asymptotically. 
They also demonstrated, using simulation experiments, 
that the test is correctly sized even in 
finite samples and possesses high powers against various alternatives.

In Remark 2.4, GHM mentioned that $\bm{V}$ 
is positive definite if $p=1$ and $h=2$; however, 
they did not prove that $\bm{V}$ is generally positive definite. 
GHM emphasized that $\bm{V}$ can be 
estimated consistently even if it is not positive definite, and the test 
can be implemented regardless of positive definiteness of $\bm{V}$.

This paper proves that $\bm{V}$ is 
positive definite in general, which guarantees that the consistently estimated $\bm{V}$ is positive definite asymptotically. 
This implies that it is legitimate, at least asymptotically, 
to apply a decomposition for positive definite matrices, 
such as the Cholesky decomposition, to the estimated $\bm{V}$. 
This may facilitate the simulation used for computing the $p$ value.  
If $\bm{V}$ is not guaranteed to be positive definite, we need to use  
a decomposition that works even for positive semi definite matrices, 
such as eigenvalue decomposition, 
which is much slower to implement than the
Cholesky decomposition.  The paper also proposes an alternative estimator for 
$\bm{V}$, which is positive definite even in finite sample size.  
This or a similar positive definiteness property have not been   
proved for the estimator of $\bm{V}$ proposed by GHM.

The paper proceeds as follows. In the next section, it is shown that 
$\bm{V}$ is positive definite for general $p$ and $h$. 
An estimator for $\bm{V}$, which is positive definite 
even in finite sample size, is also proposed. 
The last section provides concluding remarks. The Appendix 
presents several proofs.

\section{Main results\label{MR}}
Let $\bm{X}_t=[\bm{z}_t^{\prime}, \bm{x}_t^{\prime}]^{\prime}$. 
Define 
\[
\begin{array}{llll}
\bm{\Gamma}=E(\bm{X}_t\bm{X}_t^{\prime}),&
\bm{\Gamma}_{zz}=E(\bm{z}_t \bm{z}_t^{\prime}), &
\bm{\Gamma}_{xx}=E(\bm{x}_t \bm{x}_t^{\prime}), &
\bm{\Gamma}_{zx}=E(\bm{z}_t \bm{x}_t^{\prime}),
\\
\bm{\Lambda}=E(\epsilon_t^2 \bm{X}_t\bm{X}_t^{\prime}), &
\bm{\Lambda}_{zz}=E(\epsilon_t^2 \bm{z}_t \bm{z}_t^{\prime}), & 
\bm{\Lambda}_{xx}=E(\epsilon_t^2 \bm{x}_t \bm{x}_t^{\prime}),&
\bm{\Lambda}_{zx}=E(\epsilon_t^2 \bm{z}_t \bm{x}_t^{\prime}), 
\\
\gamma_{ij}=E(x_{it}x_{jt}),& 
\bm{\gamma}_{iz}=E(x_{it}\bm{z}_t),& 
\lambda_{ij}= E(\epsilon_t^2 x_{it}x_{jt}),&
\bm{\lambda}_{iz}=E(\epsilon_t^2 x_{it}\bm{z}_t). 
\end{array}
\]
All expected values defined here and hereafter are assumed to exist and be finite. 
We also make the following assumption throughout the paper.
\\ \\
\textbf{Assumption 1}. Assumptions 2.1---2.3 in GHM hold.
\\ \\
The above assumption ensures that $\bm{\Gamma}$ is positive definite. 

GHM showed that the asymptotic covariance matrix of 
$\sqrt{n}\widehat{\bm{\beta}}_n$, namely, $\bm{V}$, is expressed 
under the null hypothesis as:
\begin{equation}
\label{V}
\bm{V} = \bm{R}\bm{S}\bm{R}^{\prime},
\end{equation}
where
\[
\bm{R}=\left[\begin{array}{ccccccc}
\bm{0}_{1\times p} & 1 & \bm{0}_{1\times p} & 0 & \cdots & 
\bm{0}_{1\times p} & 0
\\
\bm{0}_{1\times p} & 0 & \bm{0}_{1 \times p} & 1 & \cdots & 
\bm{0}_{1\times p} & 0
\\
 \vdots & \vdots & \vdots & \vdots \ddots & \vdots & \vdots
\\
\bm{0}_{1\times p} & 0 & \bm{0}_{1\times p} & 0 & \cdots & 
\bm{0}_{1\times p} & 1
\end{array}\right],
\quad
\bm{S}=\left[\begin{array}{ccc} 
\bm{\Sigma}_{11} & \cdots & \bm{\Sigma}_{1h}
\\
\vdots & \ddots & \vdots
\\
\bm{\Sigma}_{h1} & \cdots & \bm{\Sigma}_{hh}
\end{array} \right],
\]
\[
\begin{array}{l}
\bm{\Sigma}_{ij}=\bm{\Gamma}_{ii}^{-1}
\bm{\Lambda}_{ij} \bm{\Gamma}_{jj}^{-1},\quad
\bm{\Gamma}_{ij}=E[\bm{X}_{it}\bm{X}_{jt}^{\prime}],\quad
\bm{\Lambda}_{ij}=E[\epsilon_t^2 \bm{X}_{it} \bm{X}_{jt}^{\prime}], \quad
\bm{X}_{it}=[\bm{z}_t^{\prime}, x_{it}]^{\prime}, 
\quad i,j \in {1, ..., h}.
\end{array}
\]
From this expression, it is unclear whether $\bm{V}$ is generally positive definite.\footnote{$\bm{S}$ is positive semidefinite for general $p$ and $h$. See Remark 2.3 in GHM.} 
In Remark 2.4, GHM mentioned that they proved that $\bm{V}$ is 
positive definite only in the case that $p=1$ and $h=2$. The following proposition shows 
that $\bm{V}$ is generally positive definite for any $p$ and $h$.  

\begin{proposition}
\label{P1} 
The asymptotic covariance matrix $\bm{V}$ is expressed as:
\begin{equation}
\label{Vexp}
\bm{V} =
\bm{D}^{-1}( \bm{\Omega}_{zx}^{\prime}
\bm{\Omega}_{zz}\bm{\Omega}_{zx}
+ \bm{\Lambda}_{xx}- \bm{\Lambda}_{zx}^{\prime}
\bm{\Lambda}_{zz}^{-1}\bm{\Lambda}_{zx}
 ) 
\bm{D}^{-1},
\end{equation}
where $\bm{D}$ is an $h\times h$ diagonal matrix whose 
$i$th diagonal element $d_{ii}$  
is
$d_{ii}\equiv \gamma_{ii} -
\bm{\gamma}_{iz}^{\prime}\bm{\Gamma}_{zz}^{-1}
\bm{\gamma}_{iz}$, and matrices 
$\bm{\Omega}_{xz}$ and $\bm{\Omega}_{zz}
$ are:
\begin{equation}
\label{Q}\bm{\Omega}_{zx}
=[\bm{\Gamma}_{zx}^{\prime}, \bm{\Lambda}_{zx}^{\prime}
]^{\prime},
\quad\mbox{and}\quad
\bm{\Omega}_{zz}
=
\left[\begin{array}{cc}
\bm{\Gamma}_{zz}^{-1}\bm{\Lambda}_{zz}
\bm{\Gamma}_{zz}^{-1}
 & - \bm{\Gamma}_{zz}^{-1}
\\
- \bm{\Gamma}_{zz}^{-1} & \bm{\Lambda}_{zz}^{-1}
\end{array} \right].
\end{equation}
Furthermore, $\bm{V}$ is positive definite if  
$\bm{\Lambda}$ is positive definite.
\end{proposition}
\textbf{Proof}. See the Appendix.
\\ \\
Note that $d_{ii}$ is equal to the conditional variance of 
$x_{it}$, conditional on $\bm{z}_t$, when $\bm{X}_{it}$ is jointly normally distributed.

An immediate corollary of Proposition \ref{P1} is as follows (cf. Remark 2.2. in GHM). 
\begin{corollary} (conditional homoscedasticity)
When $\bm{\Lambda}=\sigma^2 \bm{\Gamma}$, $\bm{V}$ reduces to:
\begin{equation}
\label{homoV}
\bm{V}=\sigma^2\bm{D}^{-1}(
\bm{\Gamma}_{xx}-\bm{\Gamma}_{zx}^{\prime}
\bm{\Gamma}_{zz}^{-1} \bm{\Gamma}_{zx}
)\bm{D}^{-1}.
\end{equation}

\end{corollary}
\textbf{Proof}. 
When $\bm{\Lambda}=\sigma^2\bm{\Gamma}$, 
matrices $\bm{\Lambda}_{zz}$, $\bm{\Lambda}_{zx}$, and $\bm{\Lambda}_{xx}$ 
are expressed as
$\bm{\Lambda}_{zz}=\sigma^2\bm{\Gamma}_{zz}$, 
$\bm{\Lambda}_{zx}=\sigma^2\bm{\Gamma}_{zx}$, and 
$\bm{\Lambda}_{xx}
=\sigma^2
\bm{\Gamma}_{xx}$, respectively, and we have: 
\begin{equation}
\label{zeroOmg}
\begin{array}{l}
\bm{\Omega}_{zx}^{\prime}
\bm{\Omega}_{zz}\bm{\Omega}_{zx}
=
\left[\begin{array}{cc} \bm{\Gamma}_{zx}^{\prime} & 
\sigma^2\bm{\Gamma}_{zx}^{\prime} \end{array}\right]
\left[\begin{array}{cc} 
\sigma^2\bm{\Gamma}_{zz}^{-1} & -\bm{\Gamma}_{zz}^{-1}
\\
-\bm{\Gamma}_{zz}^{-1} & 
\sigma^{-2} \bm{\Gamma}_{zz}^{-1}
\end{array}\right]
\left[\begin{array}{c} \bm{\Gamma}_{zx}
\\ \sigma^2\bm{\Gamma}_{zx}
\end{array} \right]
\\
\qquad\qquad\quad
=\left[\begin{array}{cc}
\bm{0}_{h\times p} & \bm{0}_{h\times p} \end{array} \right]
\left[\begin{array}{c} \bm{\Gamma}_{zx}
\\ \sigma^2\bm{\Gamma}_{zx}
\end{array} \right]
\\
\qquad\qquad\quad
=
\bm{0}_{h\times h}.
\end{array}
\end{equation}
From (\ref{Vexp}) and (\ref{zeroOmg}), we obtain the result in (\ref{homoV}). $\Box$
\\ \\
The expression of $\bm{V}$ in (\ref{homoV}) 
is, in fact, different from the asymptotic covariance matrix of 
$\sqrt{n}\widehat{\bm{b}}_n$ in the case of conditional 
homoscedasticity, namely, 
$\sigma^2(
\bm{\Gamma}_{xx}-\bm{\Gamma}_{zx}^{\prime}
\bm{\Gamma}_{zz}^{-1} \bm{\Gamma}_{zx}
)^{-1}$, 
where $\widehat{\bm{b}}_n$ is the OLS estimator of 
$\bm{b}$ in the full regression model in (\ref{full}). 

Next, we consider the estimation of $\bm{V}$. Define $\bm{\theta}_i = [\bm{\alpha}_i^{\prime}, \beta_i]^{\prime}$.   
Let $\widehat{\bm{\theta}}_{ni}$ denote the OLS estimator for $\bm{\theta}_i$ 
in the $i$th parsimonious regression: 
$y_i = \bm{X}_{it}^{\prime}\bm{\theta}_i + u_{it}$. 
To estimate $\bm{V}$, GHM proposed the following estimator:
\[
\widehat{\bm{V}}_n=\bm{R}\widehat{\bm{S}}\bm{R}^{\prime},
\]
where $\widehat{\bm{S}}=[\widehat{\bm{\Sigma}}_{ij}]_{i,j}$,  
$\widehat{\bm{\Sigma}}_{ij}=\widehat{\bm{\Gamma}}_{ii}^{-1}
\widehat{\bm{\Lambda}}_{ij}\widehat{\bm{\Gamma}}_{jj}^{-1}$, 
$\widehat{\bm{\Gamma}}_{ij}=n^{-1}\sum_{t=1}^n\bm{X}_{it}\bm{X}_{jt}^{\prime}$, 
$\widehat{\bm{\Lambda}}_{ij}=n^{-1}\sum_{t=1}^n\widehat{u}_{it}^2
\bm{X}_{it}\bm{X}_{jt}^{\prime}$, and 
$\widehat{u}_{it}=y_t-\bm{X}_{it}^{\prime}\widehat{\bm{\theta}}_{ni}$. 
Under Assumption 1, $\widehat{\bm{V}}_n$ is a consistent estimator 
for $\bm{V}$ under $H_0$, and converges in probability to a matrix  
under $H_1$ (Theorem 2.3 in GHM).
Because it is a consistent estimator for $\bm{V}$ under $H_0$ 
and $\bm{V}$ is positive definite, as shown in Proposition \ref{P1} above,  
$\widehat{\bm{V}}_n$ is asymptotically positive definite under $H_0$.  
However, the positive definiteness may not be guaranteed in 
finite sample because 
it uses different estimates of the error term in constructing 
$\widehat{\bm{\Lambda}}_{ij}$, 
depending on $i$. 

This paper proposes a slightly different consistent 
estimator for $\bm{V}$ as: 
\[
\widetilde{\bm{V}}_n=\bm{R}\widetilde{\bm{S}}\bm{R}^{\prime},
\]  
where $\widetilde{\bm{S}}=[\widetilde{\bm{\Sigma}}_{ij}]_{i,j}$, 
$\widetilde{\bm{\Sigma}}_{ij}
=\widehat{\bm{\Gamma}}_{ii}^{-1}\widetilde{\bm{\Lambda}}_{ij}
\widehat{\bm{\Gamma}}_{jj}^{-1}$, 
$\widetilde{\bm{\Lambda}}_{ij}=n^{-1}\sum_{t=1}^n \widetilde{u}_t^2 
\bm{X}_{it}\bm{X}_{jt}^{\prime}$,  
$\widetilde{u}_t=y_t-\bm{z}_t\widetilde{\bm{\alpha}}_n$, 
and $\widetilde{\bm{\alpha}}_n$ is 
the OLS estimator of $\bm{\alpha}$ in the regression:  
$y_t = \bm{z}_t^{\prime}\bm{\alpha}+ \epsilon_t$, 
which is the regression obtained under $H_0$.
The estimator $\widetilde{\bm{V}}_n$ is almost always 
positive definite in practice, even in finite sample size. 
To see this, note that $\widetilde{u}_t$, which is the estimate 
of the error term $\epsilon_t$ under $H_0$, does not depend on $i$, and  
is common in all $\widetilde{\bm{\Lambda}}_{ij}$. 
Then, it follows from the proof of Proposition \ref{P1} that 
$\widetilde{\bm{V}}_n$ can be expressed as:
\[
\widetilde{\bm{V}}_n=\widehat{\bm{D}}^{-1}
( \widetilde{\bm{\Omega}}_{zx}^{\prime}
\widetilde{\bm{\Omega}}_{zz}
\widetilde{\bm{\Omega}}_{zx}
+ \widetilde{\bm{\Lambda}}_{xx} 
- \widetilde{\bm{\Lambda}}_{zx}^{\prime}
\widetilde{\bm{\Lambda}}_{zz}^{-1} 
\widetilde{\bm{\Lambda}}_{zx}
 ) 
\widehat{\bm{D}}_n^{-1},
\]
where 
$\widetilde{\bm{\Lambda}}_{xx}=n^{-1}\sum_{t=1}^n
\widetilde{u}_t^2 \bm{x}_t \bm{x}_t^{\prime}$, 
$\widetilde{\bm{\Lambda}}_{zx}=n^{-1}\sum_{t=1}^n
\widetilde{u}_t^2 \bm{z}_t \bm{x}_t^{\prime}$, 
$\widetilde{\bm{\Lambda}}_{zz}=n^{-1}\sum_{t=1}^n 
\widetilde{u}_t^2 \bm{z}_t \bm{z}_t^{\prime}$, 
and 
$\widehat{\bm{D}}$, $\widetilde{\bm{\Omega}}_{zx}$, 
$\widetilde{\bm{\Omega}}_{zz}$, 
are defined similarly to 
$\bm{D}$, $\bm{\Omega}_{zx}$, $\bm{\Omega}_{zz}$, respectively,   
where $\bm{\Lambda}_{zx}$ and $\bm{\Lambda}_{zz}$  are replaced with  
$\widetilde{\bm{\Lambda}}_{zx}$, and
$\widetilde{\bm{\Lambda}}_{zz}$, respectively, and
$\bm{\Gamma}_{zx}$, 
$\bm{\Gamma}_{zz}$, $\gamma_{ii}$, and $\gamma_{iz}$ are 
replaced with their corresponding sample moments, 
$\widehat{\bm{\Gamma}}_{zx}$, 
$\widehat{\bm{\Gamma}}_{zz}$, 
$\widehat{\gamma}_{ii}$, and $\widehat{\gamma}_{iz}$; here, for example, 
$\widehat{\bm{\Gamma}}_{zx}$ is defined as 
$n^{-1}\sum_{t=1}^n \bm{z}_t \bm{x}_t^{\prime}$. 
The proof of Proposition \ref{P1} implies that 
$\widetilde{\bm{V}}_n$ is positive definite whenever 
$\widetilde{\bm{\Lambda}}
=n^{-1}\sum_{t=1}^n \widetilde{u}_t^2 \bm{X}_t\bm{X}_t^{\prime}$
is positive definite. This holds true, for example, when 
$\widetilde{u}_t\not=0$ for all $t=1,...,n$,
which is practically always satisfied.\footnote{
The matrix $\widetilde{\bm{\Lambda}}$ is expressed as 
$\widetilde{\bm{\Lambda}}=n^{-1}\bm{X}^{\prime}\bm{D}_u
(\bm{X}^{\prime}\bm{D}_u)^{\prime}$, where 
$\bm{X}\equiv [\bm{X}_1$, $\cdots$, $\bm{X}_n]^{\prime}$ and 
$\bm{D}_u$ is the $n\times n$ diagonal matrix whose $i$th diagonal element 
is equal to $\widetilde{u}_t$. 
When $\widetilde{u}_t\not=0$ for all $t=1,...,n$, $\bm{D}_u$ is full rank 
and thus $\mathrm{rank}(\widetilde{\bm{\Lambda}})
=\mathrm{rank}( \bm{X}^{\prime}\bm{D}_u)=\mathrm{rank}(\bm{X}^{\prime})$. 
Then, the positive definiteness of $\widetilde{\bm{\Lambda}}$
follows from Assumption 2.2. in GHM, namely, $\bm{X}$ is of full column rank $p+h$  
almost surely.  }
\footnote{Actually,  
up to $n-(p+h)$ values of $\widetilde{u}_t$ among the $n$ values of $\widetilde{u}_t$, $t=1,...,n$,  
can be zero, depending on $\bm{X}$. For example, any  
$n-(p+h)$ diagonal elements of $\bm{D}_u$ can be zero as long as 
$p+h$ column vectors of $\bm{X}^{\prime}\bm{D}_u$
are linearly independent.}
Moreover, we obtain the following proposition on the consistency property of 
$\widetilde{\bm{V}}_n$, which is similar to Theorem 2.3 in GHM. 

\begin{proposition}
\label{P2} 
Under Assumption 1, it holds that, as $n \rightarrow \infty$, 
$\widetilde{\bm{V}}_n\overset{p}{\rightarrow} \bm{V}$ under $H_0$, and  
$\widetilde{\bm{V}}_n \overset{p}{\rightarrow} \bar{\bm{V}}$ under $H_1$, 
where $\bar{\bm{V}}$ is a matrix of positive definite. 
\end{proposition}
\textbf{Proof}. The proof is similar to the proof for Theorem 2.3 in GHM, and is thus omitted here. 
\\ 

One may concern that the performance of the max test can be affected by 
the choice of the estimator for $\bm{V}$. 
Therefore, we conducted a simulation experiment to check whether there are 
any differences in the performances of the max tests with 
$\widetilde{\bm{V}}_n$ and $\widehat{\bm{V}}_n$. 
Due to space considerations, 
we do not report detailed simulation results, but 
the unreported simulation experiment confirmed 
that the performance of the max test  
is not affected when replacing $\widehat{\bm{V}}_n$ with 
$\widetilde{\bm{V}}_n$ in terms of both size and power. 
Using an estimator for $\bm{V}$ that is positive definite  
even for finite sample size may facilitate the computation of 
the simulated $p$ value. For example, we can use the estimate $\widetilde{\bm{V}}_n$
when the estimate $\widehat{\bm{V}}_n$ is not positive definite.

\section{Concluding remarks}
This paper showed that the asymptotic covariance matrix of 
the OLS estimators in the parsimonious regressions considered in  
\citet{GHM01} is generally positive definite. 
This result may facilitate the computation of the simulated $p$ value 
needed to implement the max test. 
One may consider applying a bootstrap method to the max test. 
However, to achieve an asymptotic refinement by bootstrap, 
it is crucial that the test statistic is (asymptotically) pivotal 
(see \citet{HM01} on the asymptotic validity of a bootstrap method for the max test 
in a different context). 
Although the max test statistic in \citet{GHM01} is not even asymptotically pivotal, 
the results in this paper may be utilized to construct a pivotal 
version of the max test.  
%\\
%\\
%\\
\newpage
\Large\textbf{Appendix: Proofs }\normalsize
\\
\\
First, we state the following Lemma, which is used 
in the proof of Proposition \ref{P1} below.
\begin{lemma}
\label{L1}
If $\bm{A}$ is positive definite, then 
$\left[ \begin{array}{cc} \bm{B}^{\prime}\bm{A}\bm{B} & 
\bm{B}^{\prime} \\ \bm{B} & \bm{A}^{-1} \end{array} \right]$
is positive semidefinite.
\end{lemma}
\textbf{Proof}. For any vector $\bm{x}=[\bm{x}_1^{\prime}, 
\bm{x}_2^{\prime}]^{\prime}$, we have
\[
\begin{array}{l}
[\bm{x}_1^{\prime}, \bm{x}_2^{\prime}]
\left[ \begin{array}{cc} \bm{B}^{\prime}\bm{A}\bm{B} & 
\bm{B}^{\prime} \\ \bm{B} & \bm{A}^{-1} \end{array} \right]
\left[\begin{array}{c} \bm{x}_1 \\ \bm{x}_2 \end{array}\right]
\\
=\bm{x}_1^{\prime} 
\bm{B}^{\prime}\bm{A}\bm{B} \bm{x}_1
+\bm{x}_1^{\prime} \bm{B}^{\prime} \bm{x}_2
+\bm{x}_2^{\prime} \bm{B} \bm{x}_1
+\bm{x}_2^{\prime} \bm{A}^{-1} \bm{x}_2
\\
=(\bm{x}_1^{\prime}\bm{B}^{\prime}\bm{A}+\bm{x}_2^{\prime})
(\bm{B}\bm{x}_1+\bm{A}^{-1}\bm{x}_2)
\\
=(\bm{x}_1^{\prime}\bm{B}^{\prime}+\bm{x}_2^{\prime}\bm{A}^{-1})
\bm{A}(\bm{B}\bm{x}_1+\bm{A}^{-1}\bm{x}_2)
\\
=(\bm{B}\bm{x}_1+\bm{A}^{-1}\bm{x}_2)^{\prime}
\bm{A}(\bm{B}\bm{x}_1+\bm{A}^{-1}\bm{x}_2)\geq 0.
\end{array}
\]
Because $\bm{A}$ is positive definite, the above quadratic form  is zero if and only if 
$\bm{B}\bm{x}_1+\bm{A}^{-1}\bm{x}_2=\bm{0}$. Thus 
$\bm{x}^{\prime} \left[ \begin{array}{cc} \bm{B}^{\prime}\bm{A}\bm{B} & 
\bm{B}^{\prime} \\ \bm{B} & \bm{A}^{-1} \end{array} \right] \bm{x}=0$, when 
$\bm{x}=[\bm{x}_1^{\prime}, 
-\bm{x}_1^{\prime}\bm{B}^{\prime} \bm{A}]^{\prime}$, which shows that 
the matrix is positive semidefinite. $\Box$
\\ \\
\textbf{Proof of Proposition 1}.
\\
Let $\bm{g}_i$ be the ($p+1$)th row of $\bm{\Gamma}_{ii}^{-1}$.
By the matrix inversion formula, we obtain:
%\[
%\bm{\Gamma}_{ii}^{-1}
%=\left[ \begin{array}{cc}
%\bm{\Gamma}_{zz} & \bm{\gamma}_{iz} 
%\\
%\bm{\gamma}_{iz}^{\prime} & \gamma_{ii}
%\end{array}\right]^{-1}
%=\left[\begin{array}{cc}
%\bm{\Gamma}_{z|i}^{-1} & 
%- \bm{\Gamma}_{zz}^{-1}\bm{\gamma}_{iz} d_{ii}^{-1}
%\\
%-d_{ii}^{-1}\bm{\gamma}_{iz}^{\prime} \bm{\Gamma}_{zz}^{-1}
%&
%d_{ii}^{-1}
%\end{array}\right],
%\]
%where 
%$\bm{\Gamma}_{z|i}= \bm{\Gamma}_{zz}
%-\gamma_{ii}^{-1}\bm{\gamma}_{iz}\bm{\gamma}_{iz}^{\prime}$ and $d_{ii} = \gamma_{ii} -
%\bm{\gamma}_{iz}^{\prime}\bm{\Gamma}_{zz}^{-1}
%\bm{\gamma}_{iz}$. 
%Therefore, we have:
\begin{equation}
\label{invg}
\bm{g}_i =\left[\begin{array}{cc}
-d_{ii}^{-1} 
\bm{\gamma}_{iz}^{\prime}\bm{\Gamma}_{zz}^{-1},
&
d_{ii}^{-1}
\end{array}\right], 
\end{equation}
and $\bm{V}$ is rewritten using $\bm{g}_i$ as:
\begin{equation}
\begin{array}{l}
\label{VgAg}
\bm{V}
=
\bm{R}\bm{S}\bm{R}^{\prime}
=\left[\begin{array}{cccc}
\bm{g}_{1}\bm{\Lambda}_{11}\bm{g}_{1}^{\prime} 
& \bm{g}_{1}\bm{\Lambda}_{12}\bm{g}_{2}^{\prime} & \cdots & 
\bm{g}_{1}\bm{\Lambda}_{1h}\bm{g}_{h}^{\prime} 
\\
\bm{g}_{2}\bm{\Lambda}_{21}\bm{g}_{1}^{\prime} & 
\bm{g}_{2}\bm{\Lambda}_{22}\bm{g}_{2}^{\prime} & \cdots & 
\bm{g}_{2}\bm{\Lambda}_{2h}\bm{g}_{h}^{\prime} \\
 \vdots & \vdots & \ddots & \vdots \\
\bm{g}_{h}\bm{\Lambda}_{h1}\bm{g}_{1}^{\prime} & 
\bm{g}_{h}\bm{\Lambda}_{h2}\bm{g}_{2}^{\prime} & \cdots & 
\bm{g}_{h}\bm{\Lambda}_{hh}\bm{g}_{h}^{\prime} 
 \end{array} \right].
\end{array}
\end{equation}
From (\ref{invg}) and noting that
$\bm{\Lambda}_{ij} = \left[\begin{array}{cc}
\bm{\Lambda}_{zz} & \bm{\lambda}_{jz}
\\
\bm{\lambda}_{iz}^{\prime} & \lambda_{ij}
\end{array} \right]
$, we have:
\begin{equation}
\label{gAg}
\begin{array}{l}
\bm{g}_{i}\bm{\Lambda}_{ij}
\bm{g}_{j}
\\
=\left[\begin{array}{cc}
-d_{ii}^{-1} \bm{\gamma}_{iz}^{\prime} \bm{\Gamma}_{zz}^{-1}
& d_{ii}^{-1}
\end{array}\right]
\left[\begin{array}{cc}
\bm{\Lambda}_{zz} & \bm{\lambda}_{jz}
\\ \bm{\lambda}_{iz}^{\prime} & \lambda_{ij}
\end{array} \right]
\left[ \begin{array}{c}
-d_{jj}^{-1} \bm{\Gamma}_{zz}^{-1} \bm{\gamma}_{jz}
\\ d_{jj}^{-1}
\end{array} \right]
\\
= d_{ii}^{-1} d_{jj}^{-1} \left[\begin{array}{l}
\bm{\gamma}_{iz}^{\prime} \bm{\Gamma}_{zz}^{-1}
\bm{\Lambda}_{zz}\bm{\Gamma}_{zz}^{-1}\bm{\gamma}_{jz} 
-\bm{\lambda}_{iz}^{\prime}\bm{\Gamma}_{zz}^{-1}\bm{\gamma}_{jz}
-\bm{\gamma}_{iz}^{\prime} \bm{\Gamma}_{zz}^{-1}
\bm{\lambda}_{jz} +\lambda_{ij}
\end{array}
\right] 
\\
=d_{ii}^{-1} d_{jj}^{-1} \left\{ 
\left[ \begin{array}{cc} \bm{\gamma}_{iz}^{\prime} &
\bm{\lambda}_{iz}^{\prime} \end{array} \right]
\left[\begin{array}{cc}
\bm{\Gamma}_{zz}^{-1}
\bm{\Lambda}_{zz}\bm{\Gamma}_{zz}^{-1} & 
- \bm{\Gamma}_{zz}^{-1}
\\
- \bm{\Gamma}_{zz}^{-1} & \bm{\Lambda}_{zz}^{-1} 
\end{array} \right]
\left[ \begin{array}{c} \bm{\gamma}_{jz} \\
\bm{\lambda}_{jz} \end{array} \right]
+\lambda_{ij}-\bm{\lambda}_{iz}^{\prime} \bm{\Lambda}_{zz}^{-1}
\bm{\lambda}_{jz}
\right\}
\\
=d_{ii}^{-1} d_{jj}^{-1} \left\{ 
\bm{\omega}_{iz}^{\prime}\bm{\Omega}_{zz}\bm{\omega}_{iz}
+\lambda_{ij} -\bm{\lambda}_{iz}^{\prime} \bm{\Lambda}_{zz}^{-1}
\bm{\lambda}_{jz}
\right\},
\end{array}
\end{equation}
where $\bm{\omega}_{iz}=\left[ \begin{array}{cc} \bm{\gamma}_{iz}^{\prime} &
\bm{\lambda}_{iz}^{\prime} \end{array} \right]^{\prime}$.  
From (\ref{VgAg}) and (\ref{gAg}), we obtain:
\[
\begin{array}{l}
\bm{V}=\bm{D}^{-1}
\left\{ 
\left[\begin{array}{cccc}
\bm{\omega}_{1z}^{\prime}\bm{\Omega}_{zz}\bm{\omega}_{1z}
& \cdots & \bm{\omega}_{1z}^{\prime}\bm{\Omega}_{zz}\bm{\omega}_{hz}
\\
\vdots & \ddots & \vdots
\\
\bm{\omega}_{hz}^{\prime}\bm{\Omega}_{zz}\bm{\omega}_{1z}
& \cdots & \bm{\omega}_{hz}^{\prime}\bm{\Omega}_{zz}\bm{\omega}_{hz}
\end{array}\right] 
+\left[ \begin{array}{ccc}
\lambda_{11} & \cdots & \lambda_{1h} \\ 
\vdots & \ddots & \vdots \\
\lambda_{h1} & \cdots & \lambda_{hh}
\end{array}\right] 
\right.
\\
\qquad\qquad\qquad \left.
-
\left[\begin{array}{cccc}
\bm{\lambda}_{1z}^{\prime} \bm{\Lambda}_{zz}^{-1}\bm{\lambda}_{1z}
& \cdots & \bm{\lambda}_{1z}^{\prime}\bm{\Lambda}_{zz}^{-1}
\bm{\lambda}_{hz}
\\ \vdots & \ddots & \vdots
\\ \bm{\lambda}_{hz}^{\prime}\bm{\Lambda}_{zz}^{-1}\bm{\lambda}_{1z}
& \cdots & \bm{\lambda}_{hz}^{\prime}\bm{\Lambda}_{zz}^{-1}\bm{\lambda}_{hz}
\end{array}\right]
\right\}\bm{D}^{-1}
\\
\quad
=\bm{D}^{-1}
\left\{ 
\left[\begin{array}{c}
\bm{\omega}_{1z}^{\prime}\\ \vdots \\ \bm{\omega}_{hz}^{\prime}
\end{array}\right] \bm{\Omega}_{zz}
\left[\begin{array}{ccc}
\bm{\omega}_{1z}&\cdots&\bm{\omega}_{hz}
\end{array}\right]
+\bm{\Lambda}_{xx}
- \left[\begin{array}{c} \bm{\lambda}_{1z}^{\prime} \\ \vdots
\\ \bm{\lambda}_{hz}^{\prime} \end{array}\right]
\bm{\Lambda}_{zz}^{-1}
\left[\begin{array}{ccc} \bm{\lambda}_{1z} & \cdots &
\bm{\lambda}_{hz} \end{array}\right]
\right\}\bm{D}^{-1}
\\
\quad =\bm{D}^{-1} 
(\bm{\Omega}_{zx}^{\prime}\bm{\Omega}_{zz}\bm{\Omega}_{zx} 
+ \bm{\Lambda}_{xx} - \bm{\Lambda}_{zx}^{\prime}
\bm{\Lambda}_{zz}^{-1} \bm{\Lambda}_{zx} )
 \bm{D}^{-1}, 
\end{array}
\]
which completes the proof of the first part of Proposition \ref{P1}. 
Lemma \ref{L1} implies that, if $\bm{\Lambda}$ is positive definite, 
then $\bm{\Omega}_{zz}$ is positive semidefinite. 
Similarly, matrix $\bm{\Lambda}_{xx}
-\bm{\Lambda}_{zx}^{\prime}\bm{\Lambda}_{zz}^{-1} \bm{\Lambda}_{zx} $
is positive definite when $\bm{\Lambda}$ is positive definite 
because it is the Schur complement of $\bm{\Lambda}$ with respect to 
$\bm{\Lambda}_{zz}$. 
These results imply that $\bm{V}$ is positive definite. $\Box$

\end{document}